\documentclass[11pt]{amsart}

\usepackage[T1]{fontenc}
\usepackage[utf8]{inputenc}
\usepackage{microtype}
\usepackage[hidelinks]{hyperref}
\usepackage{xurl}

\newcommand{\doi}[1]{\href{https://doi.org/#1}{doi:#1}}

\title[Training Mathematicians in the Age of AI]{Training Mathematicians in the Age of AI:\\
Intellectual Agency, Cognitive Offloading, and the PhD Thesis}
\author{Thomas Koberda}

\begin{document}

\begin{abstract}
Powerful artificial intelligence is weakening the traditional relationship between mathematical output and evidence of mathematical expertise. In particular, the production of an original theorem or a polished dissertation can no longer, by itself, certify the intellectual formation of its nominal author. I argue that graduate mathematical education should therefore be organized around the formation of intellectual agency: internal technical competence, mathematical judgment, understanding, and responsible participation in a shared intellectual culture. I distinguish productive from premature cognitive offloading, propose complementary independent and AI-augmented modes of training, and suggest a corresponding reformulation of the role of the PhD thesis and dissertation defense. More broadly, I argue that academic mathematics should understand itself increasingly as an institution for the reproduction and stewardship of human mathematical expertise rather than primarily as a mechanism for producing theorems.
\end{abstract}


\maketitle

With the rapid development of powerful AI capable of solving difficult open problems in mathematics and automating a substantial part of human mathematical cognition, it is incumbent on the mathematical community to articulate a vision for the future of all aspects of the profession. One of the more urgent issues is the reform and adaptation of graduate student training, and the role of the PhD thesis in a world where a substantial amount of theorem discovery and proving can be automated. In this essay, I will attempt to formulate and justify a vision that should remain robust under unforeseen and unforeseeable advances in technology.

The role of the academy in the age of AI is a broad subject with many urgent questions, and so I will limit my scope significantly. I will not address other equally important questions such as: Should we encourage students to use AI? Should we use AI ourselves to perform research, and if so, to what degree? How does one use AI ethically, and how does one teach AI ethics? Can AI participate meaningfully in human mathematical life? How do we address inequities that will arise from, or be exaggerated by, AI? How do we convince university administrators and society more generally that research that is not entirely AI-based is still worthwhile? My aim is for the substance of my argument to be largely independent of the answers to the first five questions. For the sixth, I will justify a related set of assumptions below, particularly in axioms four and six. My hope is that the arguments here can also contribute to a broader case for the continued value of human mathematicians and academic mathematics.

A contemporaneous paper by David Glickenstein addresses many of the same immediate questions from a complementary direction \cite{Glickenstein}. Glickenstein places \emph{mathematical judgment} at the center of graduate training and develops concrete proposals for courses, qualifying examinations, oral examinations, dissertations, and defenses. There is substantial overlap with the present essay, especially in the emphasis on verification, independent command, responsible tool use, and assessment that goes beyond polished output. The emphases are nevertheless different. Glickenstein's paper is principally concerned with program design and assessment, whereas my aim is to give a broader account of what graduate education should preserve and reproduce: intellectual agency, the management of cognitive offloading, and stewardship of a living mathematical culture. I will note some more specific points of contact below.

\section{Axiomatics}

I will begin with several axiomatic points that I will treat as the basis for this essay. Certainly, some of them may be debatable, but I will attempt to justify them and to make the assumptions themselves as minimal as possible. Several echo points in the Leiden Declaration on Artificial Intelligence and Mathematics \cite{Leiden}, and can be viewed as emerging community principles for interaction with artificial intelligence. In particular, axioms one, three, and five below overlap substantially with principles articulated there.

First, \emph{a society without human experts in areas of human knowledge is undesirable}. This point has been articulated by many other authors and thinkers, and seems nearly self-justifying insofar as robustness of human society is regarded as desirable. Many analogies, nearly all imperfect, may be applied here. An autopilot can safely and reliably perform much of the work of flying a modern aircraft, yet commercial aviation continues to depend on highly trained pilots capable of supervising automated systems and intervening when circumstances require it. Engineers still carry out design, even though modern engineers do not ordinarily draw professional plans by hand or use slide rules.

If artificial intelligence is to remain answerable to human interests, humans must retain enough independent expertise to evaluate its outputs, interrogate its reasoning where possible, and govern the circumstances in which it is trusted. Human exceptionalism, and nostalgia for a time when human control over society's destiny was more secure, are not the point. Expertise provides epistemic redundancy, allowing judgments made by machines to be interrogated, evaluated, and even overruled by experts. Robustness, however, is not the only value. A society whose members can receive correct answers to questions, but cannot independently criticize, interrogate, or understand the systems on which those answers depend, has surrendered a core part of its epistemic self-government. The expert is therefore not merely a highly trained technician, but a steward of knowledge, culture, history, and values: someone who can exercise judgment, react dynamically to unusual circumstances, spot mistakes, and transmit knowledge and experience to the next generation. Experts are, in short, custodians of epistemic self-government.

Second, \emph{mere access to tools does not guarantee successful use of those tools}. In particular, the mere presence of artificial intelligence as a tool in the world will not by itself produce human mathematical understanding. A subtly incorrect argument produced by an AI may be easy to elicit while requiring years of mathematical training to diagnose. Powerful artificial intelligence therefore creates an asymmetry between generation and evaluation: eliciting sophisticated mathematical work may require little expertise, while detecting a defect in that work may require a great deal.

Tying this to the first axiom, even the expertise required to use a tool well is deliberately formed through education; nor does this axiom obviously become irrelevant with the arrival of an ``automated researcher'' or artificial general intelligence, should such things arise. If humans are to remain epistemically responsible, then humans require intelligible grounds for deciding when to accept a machine output. Machine-checkable proofs, including proofs formally verified in systems such as Lean, are not substitutes for understanding. Formal verification answers a different question from understanding: whether a formal statement follows from specified assumptions under the rules of the system, rather than why an argument works, what its essential ideas are, or what the result means.

Another way to see this is to look at mathematical practice. Rarely are mathematical ideas or arguments ``one-and-done.'' It is standard mathematical practice to revisit, reformulate, reprove, simplify, conceptualize, generalize, connect, and canonize, and in general to build mathematical intuition. This process might be called \emph{digestion}, and seems compatible with Husserl's phenomenological view of mathematical understanding. I will return to Husserl's theory below.

Third, \emph{human mathematical culture is a practice that is not exhausted by a collection of propositions}. This too has been formulated so often as to be clich\'e, often along the lines of ``mathematics is a human activity.'' Perhaps the best justification for this axiom also comes from mathematical practice. Mathematicians often think about mathematics privately, but mathematical practice is characteristically social. Mathematics is communicated through seminars and conferences, course lectures and informal discussions, and, perhaps most importantly, preserved for the future in writing. The reverse of the Fields Medal bears the inscription ``\emph{Congregati ex toto orbe mathematici ob scripta insignia tribuere}'': mathematicians from the whole world, assembled together, conferring recognition because of outstanding writings \cite{IMUFields}.

Connecting back to Husserl and the second axiom, through inscription and symbolic practice, whether human or machine, mathematical meanings become sedimented. A dynamic mathematical community reactivates these meanings in new ways through rediscovery, reinterpretation, and reuse. A conceivable future world is one with an arbitrarily deep sediment of formally verified, perfectly correct mathematical propositions, beneath a surface on which almost no human mathematical life remains.

There is also an affective dimension to mathematical practice that should not be treated as incidental. Curiosity, surprise, aesthetic pleasure, frustration followed by understanding, and the excitement of shared discovery are among the internal rewards of mathematical life. A world containing unlimited correct mathematics but almost no humans capable of understanding it would therefore lose not only epistemic redundancy, but also a distinctive form of human intellectual experience.

Fourth, \emph{expertise and intellectual culture must be deliberately reproduced across generations}. This is closely related to the second axiom. Expertise is formed through an apprenticeship-like process in which one is exposed to examples, diverse viewpoints, books, conversations, and seminars; the apprentice hones problem-solving skills, learns from failed attempts, and absorbs wisdom from others who have more experience. In the process, the learner acquires the ability to judge arguments as strong or lacking and conjectures as interesting or less so, the insight to see where genuine difficulties lie, the judgment to distinguish technical from conceptual obstacles, and the maturity to distinguish natural from artificial generalizations.

Intellectual culture is the environment in which such an apprenticeship takes place. In the absence of an intellectual environment adapted to this kind of apprenticeship, the systematic reproduction of high-level human mathematical expertise becomes extraordinarily difficult. Artificial intelligence can complement a human-centered intellectual environment, but access to artificial intelligence by itself is not a substitute for that environment and does not supplant the process of expertise formation.

This axiom also connects explicitly to Husserl's later philosophy of mathematical practice, which views science as an intergenerational and intersubjective praxis. Modern philosophy of mathematical practice similarly emphasizes the centrality of mathematical agents, history, heuristic activity, and social epistemology; see, for example, Hartimo \cite{Hartimo}, Kitcher \cite{Kitcher}, Lave and Wenger \cite{LaveWenger}, and Lakatos \cite{Lakatos}.

Fifth, \emph{artificial intelligence and its influence on mathematical practice are here to stay}. There is no plausible return to a status quo ante in which mathematical work is produced and evaluated under the assumption that substantial cognitive automation is unavailable. Individual mathematicians may choose to use artificial intelligence extensively, sparingly, or not at all, and those choices should remain theirs. But even a mathematician who never uses artificial intelligence now works in a world in which it exists. This conclusion does not require AI capabilities to continue improving. The mere availability of systems already capable of substantial mathematical assistance permanently changes what can be inferred from a finished artifact.

One important and immediate consequence is that the epistemic meaning of a finished mathematical product has changed. When a new theorem or proof appears, its existence alone no longer provides reliable evidence about the character or quantity of human intellectual labor that produced it. This remains true even when artificial intelligence is not acknowledged or, indeed, was not used: the reader of the finished work generally cannot reconstruct its provenance from the artifact itself. Norms of disclosure, as articulated by the Leiden Declaration for instance \cite{Leiden}, may be desirable and useful, but they cannot recreate the former situation in which sophisticated mathematical output was ordinarily and reasonably treated as substantial evidence of sophisticated human mathematical activity. Sathi makes an analogous observation in the context of authorship: once AI becomes a live possibility, a reader may admire a polished document while becoming uncertain whether that admiration should extend to the person whose name appears on it \cite{Sathi}.

The indelible effect of artificial intelligence on mathematics need not be regarded as a corruption of mathematics or of the academy. It does, however, mean that the mere production of solved problems can no longer serve as a sufficient basis for mathematical identity, prestige, or professional evaluation. Theorems and proofs will remain important mathematical objects, and solving problems will remain an important mathematical activity. But the existence of a solution increasingly tells us less about the mathematician who presents it. If the profession wishes to recognize and reproduce human mathematical expertise, it must therefore look beyond mathematical output alone, toward understanding, judgment, agency, synthesis, stewardship, and the ability to participate intelligently in a living mathematical culture.

Sixth, \emph{the academy is one of society's principal institutions for reproducing and stewarding expertise, and its existing culture is well suited to adaptation in a world affected by artificial intelligence}. The academy, and in particular mathematics departments and their mechanisms for training doctoral students, provide an imperfect but durable model for the kinds of apprenticeships needed to reproduce and transmit mathematical expertise and intellectual culture. Among their distinguishing features are distributed expertise across many institutions and people; the preservation of subjects that are temporarily unfashionable; some insulation from immediate commercial utility; long time horizons that foster broader perspectives; institutional memory; and teaching as an integral part of learning. The relative autonomy and security of academic departments can, at their best, allow dynamic responses to technological shifts and provide faculty with the freedom to redefine what meaningful participation in mathematical life should entail.

The claim is not that the present organization of academic mathematics is optimal. Indeed, much of this essay will argue that its incentive structures require substantial revision; nor is the claim that universities are uniquely capable of reproducing mathematical expertise. New institutions may emerge, and such developments should be welcomed. The point is pragmatic and path-dependent: academia already contains costly, distributed, multigenerational infrastructure for apprenticeship and stewardship, and it would be reckless to dismantle that infrastructure on the assumption that an equivalent culture can easily be recreated later.

I will not invoke this axiom repeatedly in what follows. It serves primarily as a framing assumption for why the academy is already well positioned to train mathematicians for an AI-rich world, and why we should begin by adapting this existing infrastructure rather than assuming that it has become obsolete.

\section{The formation of the mathematician}

What is the purpose of doing a PhD in mathematics? Many mathematicians might agree with the assertion that a PhD program trains students to operate as \emph{independent researchers}. In principle, such researchers should be able to participate fully in the mathematical community with little more than a suitable typesetting program, access to a well-stocked library, and perhaps an internet connection through which to read the arXiv. Here, by participation in the mathematical community, I mean finding and solving open problems judged by the practitioner to be worthy of attention, developing and executing research programs, and documenting and disseminating the resulting work.

Until very recently, it was conventional to regard the completion of a dissertation thesis as one of the central aims of a PhD program. This is not to impute intellectual narrowness to anyone who holds or held such a belief, since the existence of a strong thesis long served as a reasonably reliable proxy for many of the capacities necessary for participation in mathematics. Suppose a thesis contributed strong, original results and that the document itself was technically serious. One could reasonably infer that its author was technically proficient; that they possessed the persistence, discipline, organization, and openness required to acquire substantial familiarity with a field; and that they had developed some degree of judgment and mathematical maturity, presumably under the guidance of an advisor. In particular, one could infer a meaningful degree of independence.

Artificial intelligence weakens this inference. Because the inference is itself an epistemic process carried out by individuals and communities, the fifth axiom implies that the dissertation can no longer, by itself, serve as a reliable proxy, much less sufficient evidence, for the formation of an independent mathematician. At the same time, artificial intelligence exposes weaknesses in the notion of the independent researcher as it is commonly understood.

Sathi's distinction between cognitive authorship and what he calls the institutional or ``administrative'' author is useful here \cite{Sathi}. The problem is not merely whether a theorem or dissertation is ``authored'' by the candidate in the institutional sense; the question is what relationship remains between the named mathematician and the mathematical cognition represented by the artifact. Sathi asks how institutions can recover evidence of cognitive authorship when polished artifacts no longer provide it. My concern here is somewhat different: not to identify domains that artificial intelligence cannot enter, but to specify the forms of human mathematical competence and agency that graduate education should continue to reproduce even if artificial intelligence eventually enters all of them.

We therefore need to redefine not only the role of the dissertation in the formation of a mathematician, but also what we mean by an independent researcher. In what follows, I will argue that a dissertation should certify internal competence, judgment, agency, understanding, communication, and stewardship. More broadly, I will argue that the independent researcher should cease to function as our central ideal and should be replaced by the \emph{intellectual agent}.

\section{Agency after AI}

Let us continue with the notion of the independent researcher. The word ``independent'' itself is misleading, and increasingly inappropriate, because it suggests that the independent researcher is one who does things without external help. This picture is implicit in the romantic image of a solitary mathematician sitting among dusty volumes in a secluded corner of a library, feverishly recording brilliant insights; it persists even in the more tempered image of the independent researcher equipped only with a library, a \LaTeX\ editor, and the arXiv. In either case, independence in this strong sense has never accurately described the lives of most mathematicians or the actual practice of mathematics.

Researchers often enter a field during their PhD studies, under the guidance of an advisor, and encounter a discipline that appears to have been carefully digested and documented in the mathematical literature, expressed in a timeless language. That impression is misleading. Definitions have themselves been digested and refined over generations. Thousands of predecessors established the theoretical foundations from which any current research paper begins. Mathematicians learn from and are influenced by collaborators, advisors, colleagues, seminar audiences, referees, and students. Independence in any absolute sense is mythological.

The practice of mathematics consists of many activities that persist in the presence of artificial intelligence, at least in its present form. A mathematician decides which research directions are interesting and which problems are worth pursuing. While working on a problem, the mathematician must recognize what they do and do not understand, formulate possible approaches and assess their feasibility, recognize when arguments are incorrect, reject bad suggestions from both machines and people, synthesize ideas, and recognize when outside expertise is needed. Even when a theorem is proved, the work is not finished: the mathematician poses follow-up problems, reformulates ideas, reflects on the relationship between new results and previous knowledge, and defends the choices that were made. These capacities form the basis of \emph{intellectual agency}.

Implicit in intellectual agency is \emph{technical proficiency}, without which meaningful mathematical agency would be impossible. A mathematician cannot exercise judgment about a field without possessing enough of that field internally to recognize its structures, difficulties, and possibilities. One cannot reliably distinguish a conceptual obstacle from a technical one, assess whether a proposed argument is plausible, recognize when an apparent generalization is artificial, or decide when an external suggestion deserves trust without substantial technical competence. Technical proficiency should not be confused here with \emph{technical virtuosity}, nor with an ability to execute every difficult argument unaided.

Intellectual agency should therefore not be understood as a managerial capacity to direct other people or machines. It is mathematical judgment exercised from within mathematics. This point becomes especially important in the presence of powerful artificial intelligence. A researcher may be able to elicit sophisticated arguments, examples, or conjectures with very little technical knowledge, but the ability to generate such output is not itself evidence of intellectual agency. Technical proficiency supplies the internal reference point against which externally generated mathematics can be interpreted, criticized, reconstructed, and incorporated into one's own understanding. The required level and character of proficiency will vary across fields and across stages of a career, but substantial internal mathematical competence is indispensable.

The relationship between \emph{intellectual agency} and Glickenstein's \emph{mathematical judgment} is close but not one of identity. Glickenstein distinguishes local content mastery from disciplinary mathematical formation and uses mathematical judgment as an organizing construct for capacities such as validating proofs, assessing definitions and examples, monitoring understanding, transferring ideas, and governing tools \cite{Glickenstein}. This distinction parallels, without coinciding with, the distinction here between technical proficiency and intellectual agency. Mathematical judgment is a central component of intellectual agency as I use the term, but intellectual agency also includes question selection, direction of a research program, synthesis, recognition of when outside expertise is needed, and responsible participation in the intellectual network described below.
This perspective also has affinities with the competency-based tradition in mathematics education developed by Niss and Højgaard \cite{NissHojgaard}; the notion of intellectual agency used here is broader, incorporating such mathematical competencies together with the direction, evaluation, and integration of mathematical activity.

The components of intellectual agency enumerated above are probably not especially controversial. Most mathematicians would recognize these activities in their own mathematical lives, as well as the importance of imparting them to students. Even in a world where large parts of theorem proving and verification are automated, and even in a world where artificial intelligences themselves possess intellectual agency, humans may still possess intellectual agency as well. Indeed, it is difficult to imagine a human expert in mathematics who participates fully in the discipline while lacking it entirely.

Most importantly, intellectual agency is independent of reliance on machines. A human researcher can rely heavily on artificial intelligence while retaining substantial intellectual agency. Conversely, a human researcher can satisfy traditional outward markers of independence, by working largely alone and producing technically strong work, while exercising relatively little agency over question selection, interpretation, synthesis, or the relationship of their work to the broader discipline.

\section{Cognitive offloading}

Because technical proficiency is a precondition for intellectual agency, graduate education must distinguish forms of offloading that preserve or deepen proficiency from those that prevent proficiency from forming.

For mathematicians trained in the pre-AI era, intellectual agency and the technical proficiency required to participate in mathematics were often developed through long hours of engagement with mathematical problems: parsing minute details, organizing lemmata, synthesizing insights with ideas from the literature, and carefully documenting one's findings. Much of this process is dull, discouraging, tedious, or even painful. None of this implies that all such friction is educationally valuable. Some forms of difficulty build technical skill and judgment; others are simply wasted time. Once intellectual agency has developed, however, and technical skill has been acquired in one field, technical competence in other fields often becomes easier to acquire, while mathematical agency can expand with less effort.

It is not yet clear whether there is an AI-augmented shortcut to achieving intellectual agency. A major impediment to its careful development in the presence of powerful artificial intelligence arises from an incentive system built around authorship of papers containing proofs of theorems. It is already possible in some cases for sophisticated research-level mathematical output to be generated with surprisingly little mathematical input from the human operator; see, for instance, Tao's discussion of recent controlled evaluations of frontier systems \cite{Tao}. The point is crucial, and echoes what has already been argued above: \emph{an incentive system in which theorem production serves as the principal currency of credit and prestige is no longer adequate and must be reformed}. This is not to say that theorems are no longer important. Rather, if mathematics is to remain a human practice and if there are to remain human experts in mathematics, then intellectual agency itself must be cultivated, transmitted across generations, and rewarded.

The development of intellectual agency is fundamentally a cognitive task, and one broad fear surrounding artificial intelligence is that human intellect will eventually serve no function in its presence. If humans are not to cede intellectual agency to machines, however, this fear is too coarse to be useful. A better question is: with the goal of developing and maintaining intellectual agency, which cognitive tasks can safely be offloaded to artificial intelligence, and which must first be sufficiently internalized? The introduction of calculators and other computational tools has already shown that some forms of cognitive offloading can safely render particular manual competencies obsolete. The problem is therefore not to preserve every existing skill, but to determine which internal capacities remain constitutive of understanding and agency.

For a mathematician who already possesses substantial intellectual agency and technical skill, more cognitive tasks can probably be safely offloaded than for a student whose expertise is still forming. Nevertheless, a mathematician must retain a sense for when a proposition has actually been proved, an intuitive grasp of the structure of arguments, the ability to reconstruct details when necessary, and the ability to detect nonsense. The mathematician must not allow command of foundational material to atrophy, must retain the patience and persistence to reason through unfamiliar problems, and must preserve enough memory of a field to remain oriented within it.

Here it is useful to distinguish \emph{productive offloading} from \emph{premature offloading}. Productive offloading delegates a task, cognitive or otherwise, to a machine in a way that removes unnecessary drudgery, expands capacity, or allows attention to be directed elsewhere without undermining the competence on which judgment depends. Premature offloading delegates a task whose performance is itself part of the process through which the relevant competence would have formed. I use these terms for a distinction suggested by the broader cognitive-offloading literature; see Risko and Gilbert \cite{RiskoGilbert}. Sathi makes a closely related point in the context of writing, arguing that composition is often part of the production and testing of thought rather than the mere transcription of pre-existing thought; replacing that activity with AI can therefore outsource some of the thinking that writing would otherwise force one to do \cite{Sathi}.

What constitutes productive or premature offloading is likely to vary substantially among individuals and across stages of a career. For the purposes of developing intellectual agency in the age of AI, however, graduate training must minimize premature offloading. Exactly how to achieve this will likely require empirical evidence and the collective experience of the mathematics community over many years.

One reasonable place to begin is by distinguishing two modes of training: an \emph{independent mode} and an \emph{augmented mode}. In independent mode, students would work without reliance on artificial intelligence when the pedagogical aim is to develop or assess capacities that AI might otherwise perform for them. This is the traditional engagement with strategies, failed approaches, difficult details, reconstruction, and sustained independent reasoning. In augmented mode, students would be allowed, and often actively encouraged, to make careful, intelligent, and perhaps even aggressive use of artificial intelligence in pursuing their research programs.

Related distinctions are already appearing in mathematics education. Yoon, Lee, Lee, and Kwon \cite{YoonEtAl}, for example, contrast ``initial GenAI-independent proving'' with proving in which students interact with GenAI from the outset. I use the broader terms \emph{independent mode} and \emph{augmented mode} for a related distinction applied here to the entire process of graduate mathematical formation and research.

Glickenstein proposes a closely related institutional distinction: AI policy should vary with the object of assessment, with timed qualifying examinations normally AI-free when the goal is to measure independent command, while explicitly AI-permitted tasks can instead assess the ability to audit, verify, correct, and defend machine-assisted work \cite{Glickenstein}. This is close in spirit to the independent/augmented distinction, though the latter is meant as a broader description of modes of formation and research rather than as a taxonomy of assessment formats.

The independent mode has a diagnostic and developmental function; it is not intended to carry normative or moral weight. I do not suggest that AI contaminates mathematics. Rather, when a graduate program seeks evidence that a particular capacity has been internalized, it must temporarily restrict whichever external aids would perform that capacity for the student. The permitted aids should depend on what is being developed or assessed. The point is not abstinence from tools, but evidence that there remains a mathematician capable of using them.

\section{Mathematics and the intellectual network}

A theorem bearing the names of one or several mathematicians arises from a large web of dependencies. The language in which the subject is discussed was developed and refined by generations of mathematicians, as were the theoretical foundations from which a particular research paper begins. Teachers and advisors provide formal and informal training; colleagues ask questions; referees and readers correct mistakes and force clarification; seminars generate criticism and new points of view; institutions provide time. Authorship has functioned as a compressed attribution mechanism for a much broader network of dependencies. I make no claim of originality for this observation; Thurston makes closely related points in \emph{On Proof and Progress in Mathematics} \cite{Thurston}. In a world with powerful artificial intelligence, however, this networked character of mathematical production becomes increasingly difficult to ignore, and it becomes especially easy to exaggerate the independence of a particular achievement.

None of this is intended to deny individual excellence. There have been, and remain, truly remarkable mathematicians. The point is rather that mathematical excellence is made possible by the vast ecosystem that constitutes intellectual culture generally and mathematical culture in particular. If we wish to maintain a world with human mathematicians who possess intellectual agency and in which such excellence remains possible, then human participants in mathematics must act as stewards of intellectual culture. Mathematicians, especially academic mathematicians, already perform these functions to varying degrees: teaching, judging, remembering, explaining, correcting, mentoring, synthesizing, organizing, and deciding what deserves attention.

Some of the knowledge involved in this stewardship is tacit and situated rather than fully encoded in the literature: which approaches have quietly failed, which ideas are currently in vogue, which conjectures experts actually regard as plausible, or which work has been unfairly neglected. This does not mean that such information will remain inaccessible to artificial intelligence. The point is that participating in the formation, criticism, and revision of these judgments is itself part of human mathematical agency.

Stewardship should not be confused with conservatism. A steward of a living intellectual culture does not merely preserve inherited ideas, but criticizes, reorganizes, discards, and replaces them; nor should a networked conception of mathematics diminish the importance of careful attribution. Quite the opposite: recognizing the collective infrastructure behind mathematics should make us more attentive to the particular contributions of individuals, including forms of intellectual and institutional labor that traditional authorship practices have often rendered invisible.

Mathematicians judge and are judged by each other continually. The question ``What did this mathematician prove?'' can no longer function as a sufficient summary of mathematical contribution. A better question is closer to: \emph{Through your intellectual agency, what have you helped the mathematical community understand?} Here again the reader may compare Thurston \cite{Thurston}. If we wish to maintain a world with human mathematicians in possession of intellectual agency, then students must be trained to function as responsible nodes in this intellectual network.

There is a danger lurking here. Criteria such as judgment, taste, agency, and stewardship are far less legible than theorem counts, and are therefore easier to manipulate. An attempt to de-emphasize theorem production could perversely intensify prestige hierarchies, patronage, or cults of personality if senior mathematicians were simply empowered to declare whose work displays ``taste.'' Any reform of evaluative criteria must therefore combine qualitative judgment with distributed evaluation, multiple forms of evidence, transparent standards where possible, careful attribution, and resistance to any single advisor, institution, or famous mathematician functioning as the arbiter of mathematical value. The aim is to de-emphasize crude output metrics without replacing them with opaque authority.

Nothing in this discussion requires curation, synthesis, judgment, or even intellectual agency to remain uniquely human capabilities. Artificial intelligences may eventually perform any or all of these functions better than humans. The point is normative rather than competitive: if human mathematical expertise, epistemic self-government, and a living human mathematical culture are goods worth preserving, then humans must continue to exercise these capacities even where machines can exercise them as well or better. A capacity need not be uniquely human for its human exercise to be valuable.

\section{The role of the dissertation}

Tao has described one aspect of the current transition as a movement from ``proof scarcity'' toward ``proof abundance'' \cite{Tao}. This metaphor is appropriate and compatible with the arguments developed here.

By now, it should be clear that the value of a dissertation, traditionally viewed as the culmination of many years of graduate study and training, can no longer lie principally in the mere existence of a proven original theorem. Instead, the dissertation must serve as evidence of mathematical formation and intellectual agency. Original results still matter, but they must be viewed as part of a larger whole.

A student upon whom a PhD is to be conferred should, through the thesis and its defense, demonstrate an ability to perform, in a combination of independent and augmented modes, the activities that together constitute intellectual agency. These include asking worthwhile questions, understanding the relevant literature, distinguishing conceptual landmarks from technical machinery, explaining why the work matters, locating the work within a broader mathematical landscape, identifying limitations and failed approaches, and articulating plausible directions for future progress. If a student has relied substantially on artificial intelligence, they should also be able to demonstrate intelligent, transparent, and ethical use of it. In independent mode, they should be able to evaluate machine-generated arguments and to reconstruct and defend the relevant mathematics without having the machine perform those capacities for them.

In one concrete model, the dissertation document and examination would have several components. First, there should ordinarily be a substantive mathematical contribution. Demonstrated mathematical propositions are by no means irrelevant in a world with powerful AI, and remain an important element of the thesis. Second, there should be a clear intellectual narrative. Many modern dissertations already contain extensive background sections providing context for the original contribution; this should become an expected part of the thesis, explaining how the problem arose, what theoretical background is relevant, where the ideas came from, and, if relevant, which approaches failed. Third, there should be a synthesis component explaining the position of the contribution within a broader landscape and articulating what other mathematicians should learn from it. Fourth, machine-assisted methods should be discussed: which tasks were delegated to machines, what was checked and how, and where human judgment entered the process. Such disclosure should be at a level appropriate to intellectual provenance, attribution, and reproducibility, rather than becoming an exercise in forensic surveillance.

Finally, the dissertation defense should carry greater evidentiary weight than it often does at present, with equivalent supervised or extended formats available where appropriate. The defense is a venue in which the student can provide evidence not only of understanding of the dissertation but also of intellectual agency. The student might present the contents of the dissertation in a manner not unlike modern defenses, but the examination following the presentation should play a more substantial role. Examiners might ask the student to alter hypotheses, explain what fails under different assumptions, identify the key conceptual points, reconstruct part of an argument, compare alternative approaches, or articulate a vision for future work.

Glickenstein likewise argues that the dissertation defense can serve as a structural safeguard by assessing ``mathematical ownership'': a student who used AI in research or exposition should be able to explain what was used, what was checked or discarded, and why the final mathematics is correct \cite{Glickenstein}. This is closely aligned with the evidentiary role assigned to the defense here.

No single oral examination should be treated as an infallible detector of understanding. Oral performance can be affected by personality, language, disability, and comfort under pressure. The defense should therefore be understood as one part of a broader body of evidence about mathematical formation rather than as an authenticity test standing alone. The importance of assessment that remains informative in the presence of AI is already the subject of a substantial and rapidly developing literature; see, for example, Nikoli\'c and Basta Nikoli\'c \cite{Nikolic}. Sathi likewise argues for institutional mechanisms that reward evidence of cognitive engagement rather than polished output, while warning that live intellectual performance can reward charisma and speed rather than depth \cite{Sathi}.

\section{The role of graduate training}

Mathematicians have historically treated the writing of a thesis as a central aim of graduate training, in part because the thesis served as a reasonable proxy for readiness to participate independently and at a high level in the mathematical community. Once both the notion of the independent researcher and the evidentiary role of the thesis are revised, we are left with a more fundamental question: what should graduate training aim to produce?

I propose that graduate training be understood as an apprenticeship toward stewardship of the discipline and the formation of intellectual agency.

Graduate education should still involve independent problem solving. Independent-mode work is a tested means of acquiring technical proficiency, which is itself a prerequisite for intellectual agency. But many other aspects of graduate training should be recognized as central components of mathematical formation rather than treated merely as peripheral service. These include research seminars, exposition, teaching, refereeing, collaborative work, AI-assisted exploration, structured mentoring of junior students, and synthetic writing. Taken together, such activities develop capacities by which students should increasingly be evaluated: their ability to increase the mathematical community's capacity to understand, use, criticize, and extend mathematics.

Even if theorem proving becomes largely automated, we need not conclude that mathematics as a discipline has become narrower or less interesting than it was in the pre-AI era. If the incentive structure of graduate study, and of the mathematical community more generally, is reoriented toward intellectual agency and contributions to human understanding, then the range of activities recognized as serious mathematical work may expand rather than contract. Novel viewpoints, excellent expositions of difficult theories, the organization of confused literatures, conceptual synthesis, formalization, and the creation of structures that make other mathematicians more effective may all deserve substantially more professional recognition than they have traditionally received.

If AI causes technical productivity to grow rapidly, human attention and the coordination of understanding will become correspondingly scarcer. Collaboration will therefore become an increasingly important component of graduate training, together with the skills that make collaboration productive: integrating human- and machine-generated knowledge, attributing credit carefully, communicating partial understanding, disagreeing constructively, and recognizing superior judgment in others.

One serious risk in the present moment of sensational claims about the powers of artificial intelligence and the alleged death knell of the academy is that talented students who, in the pre-AI era, would have pursued graduate work in mathematics will increasingly conclude that there is no meaningful future in doing so. Even within the academy, some mathematicians have already acted on this uncertainty. Tsimerman, for example, stopped taking graduate students because of worries that a conventional research project begun now could prepare a student for a mathematical career that may no longer exist in recognizable form \cite{Hartnett}.

Potential students should not be told lies nor sold nostalgic delusions about academic mathematics. The landscape is changing rapidly, and there is no going back. Students should be disabused of the notion that being a mathematician consists primarily of producing papers containing proofs of finished results. At the same time, the potential of the future is vast. Mathematical possibility is expanding; tools capable of exploring more ideas, checking more cases, and making more connections, at speeds previously unimaginable, are suddenly available. Learning can be accelerated to a degree previous generations could scarcely have imagined. Mathematicians alive at this moment have the opportunity to help determine how a human mathematical culture will navigate this new abundance.

\section{Mathematicians and the future}

The purpose of graduate education in mathematics is therefore not principally to ensure the production of new theorems. It is to ensure the continued formation of mathematicians. A mathematician in the age of artificial intelligence should possess sufficient internal technical competence to reason without immediately delegating every difficult cognitive task; sufficient intellectual agency to decide what should be delegated and how the resulting work should be judged; and sufficient understanding to explain, criticize, reorganize, and extend the mathematics with which they work. At the same time, a mathematician must be able to participate in the social life of the discipline: to teach and learn, collaborate, attribute credit, recognize the judgment of others, and contribute to the preservation and renewal of a shared mathematical culture.

This conception should also encourage a healthier view of individual mathematical achievement. Mathematicians are not isolated producers of intellectual objects. Each of us is a node in an immense network extending backward through our teachers and the literature we inherited, outward through colleagues, students, collaborators, and readers, and forward toward mathematicians who have not yet been trained. Individual excellence is real and should be recognized, but it is made possible by this network. The corresponding professional ideal is therefore not the solitary producer of theorems, but the responsible intellectual agent and steward: someone whose work enlarges the capacity of other people to understand and do mathematics.

The PhD thesis should certify entry into this role. Its purpose is not merely to demonstrate that a candidate has caused previously unknown propositions to become known. Increasingly, the existence of such propositions will tell us too little about the person presenting them. The thesis and its defense should instead provide evidence that the candidate has acquired the competence, agency, judgment, understanding, and habits of stewardship required to participate fully in mathematical life, both with and without powerful artificial assistance.

There is no reason to regard this change as a diminishment of mathematics. If theorem production becomes abundant, the range of mathematical work that we can recognize and reward may become broader: synthesis, explanation, conceptual reorganization, teaching, criticism, collaboration, and the careful integration of new knowledge into what humans can actually understand. Artificial intelligence may radically change what mathematicians do. Graduate education should ensure that it does not eliminate the mathematician from mathematics.

\section*{Acknowledgments}

I thank Angelo Mao for thoughtful comments on an earlier draft.

\medskip
\noindent\textbf{AI disclosure.} The author made extensive use of OpenAI's ChatGPT as an interlocutor and editorial assistant while developing and revising this essay; ChatGPT provided adversarial readings and help developing and sharpening many arguments. The arguments and judgments expressed here, and the final text, are solely the author's responsibility.

\end{document}